\documentclass[12pt]{amsart}
\usepackage{amsmath,amsthm,amsfonts,amssymb}
\usepackage[all]{xy}
\usepackage{comment}

\DeclareMathOperator{\Rep}{Rep}
\DeclareMathOperator{\End}{End}
\DeclareMathOperator{\GL}{GL}
\DeclareMathOperator{\Hom}{Hom}
\DeclareMathOperator{\im}{i}

\renewcommand{\dim}{\operatorname{dim}}

\newcommand{\ve}{\varepsilon}
\newcommand{\FPdim}{{\rm FPdim}}

\newcommand{\C}{\mathbb C}
\newcommand{\Z}{\mathbb Z}
\newcommand{\Q}{\mathbb Q}

\newcommand{\CC}{\mathcal{C}}

\newcommand{\N}{\mathbb N}

\newcommand{\ot}{\otimes}

\newcommand{\B}{\mathcal{B}}

\newcommand{\om}{\omega}
\newcommand{\lan}{\langle}
\newcommand{\ra}{\rangle}
\newcommand{\la}{{\lambda}}

\newcommand{\al}{\alpha}
\newcommand{\one}{\mathbf{1}}

\newcommand{\g}{\mathfrak{g}}
\newcommand{\f}{\mathfrak{f}}
\newcommand{\e}{\mathfrak{e}}
\newcommand{\so}{\mathfrak{so}}
\newcommand{\ssl}{\mathfrak{sl}}

\newtheorem{thm}{Theorem}[section]

\newtheorem{lem}[thm]{Lemma}
\newtheorem{prop}[thm]{Proposition}

\theoremstyle{definition}
\newtheorem{defn}[thm]{Definition}

\newtheorem{conj}[thm]{Conjecture}

\begin{document}

\title[Braid representations from quantum groups]
{Braid representations from quantum groups of exceptional Lie type}
\author{Eric C. Rowell}
\address{Department of Mathematics,
    Texas A\&M University, College Station, TX 77843, USA}
   \email{rowell@math.tamu.edu}
\thanks{The author is partially supported by NSA grant H98230-08-1-0020.  Thanks go to Hans Wenzl and Matt Papanikolas for useful discussions.
}
\begin{abstract}
We study the problem of determining if the braid group representations obtained from quantum groups of types $E, F$ and $G$ at roots of unity have infinite image or not.  In particular we show that when the fusion categories associated with these quantum groups are not weakly integral, the braid group images are infinite.  This provides further evidence for a recent conjecture that weak integrality is necessary and sufficient for the braid group representations associated with any braided fusion category to have finite image.
\end{abstract}

\maketitle

\section{Introduction}
The $n$-strand braid group $\B_n$ is defined abstractly by generators $\sigma_1,\sigma_2,\ldots,\sigma_{n-1}$
satisfying relations:
\begin{enumerate}
\item[(R1)] $\sigma_i\sigma_{i+1}\sigma_i=\sigma_{i+1}\sigma_i\sigma_{i+1}$, $1\leq i\leq n-2$
\item[(R2)] $\sigma_{j}\sigma_i=\sigma_{i}\sigma_j$, $|i-j|>1$.
\end{enumerate}
Given an object $X$ in a braided fusion category $\CC$ the braiding $c_{X,X}\in\End(X^{\ot 2})$ allows one to construct a family of braid group representations via the homomorphism $\C\B_n\rightarrow \End(X^{\ot n})$ defined on the braid group generators $\sigma_i$ by
$$\sigma_i\rightarrow Id_X^{\ot i-1}\ot c_{X,X}\ot Id_X^{\ot n-i-1}.$$  For any simple object $Y\subset X^{\ot n}$ the simple $\End(X^{\ot n})$-modules $\Hom(Y,X^{\ot n})$ become $\B_n$ representations (although not necessarily irreducible).
In this paper we consider the problem of determining when the images of these representations are finite groups.
We will say a category $\CC$ has \textit{property $\mathbf{F}$} if, for all objects $X$ and all $n$, the corresponding braid group representations factor over finite groups.  For example, the representation category $\Rep(D^\om G)$ of the twisted double of a finite group $G$ always has property $\mathbf{F}$ \cite{ERW}.  Various cases related to Hecke- and BMW-algebras have been studied in the literature, see \cite{FRW,FLW,jones86,jonescmp,LRW}.  Recent interest in this question has come from the study of topological quantum computation, see \cite{FLW} for particulars.  

For any object $X$ in a fusion category one defines the \emph{FP-dimension} $\FPdim(X)$ (see \cite{ENO}) to be the largest (necessarily real) eigenvalue of the matrix representing $X$ in the left-regular representation of the Grothendieck semiring.  That is, $\FPdim(X)$ is the largest eigenvalue of the matrix whose $(i,j)$th entry is the multiplicity of the simple object $X_j$ in $X\ot X_i$.  A fusion category is said to be \emph{integral} if $\FPdim(X)\in\N$ for every object, and \emph{weakly integral} if $\FPdim(X_i)^2\in\N$ for any simple object $X_i$, or, equivalently, if $\sum_i (\FPdim(X_i))^2\in\N$ where the sum is over (representatives of isomorphism classes of) simple objects.

Empirical evidence (from the papers mentioned above, for example) partially motivates (see also \cite[Section 6]{RSW} and \cite{NR}):
\begin{conj}\label{propfconj}
 A braided fusion category $\CC$ has property $\mathbf{F}$ if, and only if, the Frobenius-Perron dimension $\FPdim(\CC)$ of $\CC$ is an integer, (\textit{i.e.} $\CC$ is \emph{weakly integral}).
\end{conj}
For example this conjecture predicts that the braid group images arising from any finite dimensional semisimple quasi-triangular quasi-Hopf algebra must be finite groups.

 The goal of this paper is to verify one direction of this conjecture for a large class of braided fusion categories--namely, we will show that if $\CC$ is a non-weakly-integral braided fusion category associated to a quantum group of exceptional Lie type at a root of unity then the braid group images are infinite.  For classical Lie types $A,B,C$ and $D$ results of this type have already appeared in \cite{FLW} and \cite{LRW}.

The body of the paper is organized into four sections.  In Section \ref{2:qgcats} we outline the necessary background for the problem.  Section \ref{3:wiclass} classifies the categories $\CC$ associated with exceptional-type quantum groups which are weakly integral and in Section \ref{4:noF} we show that property $\mathbf{F}$ fails in the non-weakly-integral cases.  We make some final remarks about the problem in Section \ref{5:conclusions}.

\section{Preliminaries}\label{2:qgcats}
Associated to any semisimple finite dimensional Lie algebra $\g$ and a complex number $q$ such that $q^2$ is a primitive $\ell$th root
of unity is a ribbon fusion category $\CC(\g,q,\ell)$.  The construction is essentially due to Andersen \cite{andersen} and his collaborators.  We refer the
reader to the survey paper \cite{Survey} and the texts \cite{BK} and \cite{Tur} for a more complete treatment.

We shall denote $q$-numbers in the standard way, i.e. $[n]:=\frac{q^n-q^{-n}}{q-q^{-1}}$.  Let $\Phi_+$ denote the
positive roots for the system of type $\g$ equipped with the form $\lan-,-\ra$ normalized so that $\lan \al,\al\ra=2$ for short roots, and $\check{\Phi}_+$ denote the positive coroots.  The dominant Weyl alcove will be denoted $P_+$ and we will denote the fundamental weights by $\la_i$. Define $m$ to be the ratio of the squared-lengths of a long root and a short root so that $m=1$ for Lie type $E$, $m=2$ for Lie type $F_4$ and $m=3$ for Lie type $G_2$.  Let $\vartheta_0$ be the highest root and $\vartheta_1$ be the highest short root.  The simple objects in $\CC(\g,q,\ell)$ are labeled by:
$$C_\ell(\g):=\begin{cases}\{\la\in P_+: \lan\la+\rho,\vartheta_0\ra<\ell\} &
\text{if $m\mid \ell$}\\
\{\la\in P_+: \lan\la+\rho,\vartheta_1\ra<\ell\} & \text{if $m\nmid
\ell$}\end{cases}$$

In this notation the $FP$-dimension for the simple object in $\CC(\g,q,\ell)$ labelled by highest weight $\lambda$ is obtained by setting
$q=e^{\pi i/\ell}$ in:
$$
\FPdim(V_{\la}):=\begin{cases} \prod_{\al\in\Phi_+}\frac{[\lan \la+\rho,\al\ra]}{[\lan \rho,\al\ra]}, & m\mid\ell\\
\prod_{\check{\al}\in\check{\Phi}_+}\frac{[\lan \la+\rho,\check{\al}\ra]}{[\lan \rho,\check{\al}\ra]}, & m\nmid\ell\end{cases}
$$
The first case is well-known, see \cite{wenzlcstar} for a proof of positivity.  The latter case is more recent, a proof is found in \cite{RJPAA}.

To deduce the infinitude of the braid group images we will employ the main result in \cite{RT} reproduced here for the reader's convenience.  Let $\varphi:\B_3\rightarrow\GL(V)$ be a $d$-dimensional irreducible representation with $2\leq d\leq 5$ and set $A=\varphi(\sigma_1)$ and $B=\varphi(\sigma_2)$.  Let $G$ denote the image  $\varphi(\B_3)=\lan A,B\ra$ i.e. the group generated by $A$ and $B$ and
define $S:=\mathrm{Spec}(A)=\mathrm{Spec}(B)$.
\begin{defn}
A linear group $\Gamma\subset\GL(V)$ is \emph{imprimitive} if $V$ is irreducible
and can be expressed as a direct sum of
subspaces $V_i$ which $\Gamma$ permutes nontrivially.  Otherwise, we say that $\Gamma$ is \emph{primitive}.
\end{defn}
\begin{thm}[\cite{RT}]\label{rttheorem}  Let $\varphi$, $G$, $A$ and $S$ be defined as
in the previous paragraph.
Let $S=\{\la_1,\ldots,\la_d\}$, and define the projective order of $A$ by
$$po(A):=\min\{t: (\la_1)^t=(\la_2)^t=\cdots=(\la_d)^t\}.$$   We use the convention
that each successive statement excludes the hypotheses of all
of the preceding cases.
\begin{enumerate}
\item[(a)] Suppose some $\la_i$ is not a root of unity, or $\la_i=\la_j$ for some $i\neq j$.  Then $G$ is \textbf{infinite}.
\item[(b)] Suppose $po(A)\leq 5$.  Then $G$ is \textbf{finite}.
\item[(c)] Suppose $G$ is imprimitive.  Then $S$ is of the form:
\begin{enumerate}
 \item[(i)] $\{\pm\chi,\alpha\}$ or $\chi\{1,\omega,\omega^2\}\cup\{\alpha\}$
with $\omega$ a primitive 3rd root of unity and $G$ is \textbf{finite} or
\item[(ii)] $\{\pm r,\pm s\}$.  In this case if $u=r/s$
is a root of unity of order $o(u)\in \{7,8,9\}\cup[11,\infty)$ then $G$ is \textbf{infinite}, if $o(u)=6$, $G$ is \textbf{finite} and if $o(u)=5$ or $10$
one cannot decide $|G|$ without further information.
\end{enumerate}
\item[(d)] Suppose $G$ is primitive. Then:
\begin{enumerate}
\item[(i)] If $d=2$ then $G$ is \textbf{infinite}.
\item[(ii)] If $d=3$ and $po(A)\geq 8$ then $G$ is \textbf{infinite}.  If $po(A)=7$, and $\frac{1}{\la_1}S$ is Galois conjugate to $\{1,e^{2\pi\im/7},e^{2k\pi\im/7}\}$ with $k$ even, $G$ is \textbf{infinite}, whereas if
$k$ is odd, $G$ is \textbf{finite}.
\item[(iii)] If $d=4$ and $po(A)\not\in\{6,\ldots,10,12,15,20,24\}$ then $G$ is \textbf{infinite}.
\item[(iv)] If $d=5$ and $po(A)\in\{7,8\}\cup[13,\infty)$ then $G$ is \textbf{infinite}.
\end{enumerate}
\end{enumerate}
\end{thm}
Notice that by (a) if $A$ is not diagonalizable then it has infinite order, so that in this case ``projective order" is a misnomer.  Since we will be using this theorem to show the image of $\B_3$ is infinite, this should not cause any confusion.

In many cases irreducibility may be deduced from:
\begin{lem}\cite[Lemma 3.2]{TW}\label{twlemma} Suppose $Z$ is a simple self-dual object such that $Z^{\ot 2}$ is isomorphic to $\bigoplus_{i=1}^d Y_i$ where $Y_i$ are distinct simple objects and such that $\sigma_1$ acts on each $Y_i$ by distinct scalars.  Then $\B_3$ acts irreducibly on the $d$-dimensional vector space $\Hom(Z,Z^{\ot 3})$.
\end{lem}
In order to use these results we must have the eigenvalues of the image of $\sigma_1$ at our disposal in order to compute the projective order of the image of $\sigma_1$.  These can be deduced from the quantum group via Reshetikhin's formula (see \cite{LeducRam}[Corollary 2.22(3)]):
Suppose that $V$ is an irreducible highest weight representation of $U_q\g$ labeled by $\la$ and $W$ is a subrepresentation of $V\ot V$ labeled by $\mu$.  Assume further that $V^{\ot 2}$ is multiplicity free.
\begin{equation}\label{eigs}
(c_{V,V})_{|_W}=\pm f(\la) q^{\lan\mu,\mu+2\rho\ra/2}\mathbf{1}_W
\end{equation}
where $f(\la)$ is an overall scale factor that depends only on $\la$ and the sign is $+1$ if $W$ appears in the symmetrization of $V\ot V$ and $-1$ if $W$ appears in the antisymmetrization of $V\ot V$.  Since the ribbon structure of $\CC(\g,q,\ell)$ is inherited from that of $Rep(U_q\g)$ we may use eqn. (\ref{eigs}) to compute the eigenvalues of the action of the generators of $\B_3$ on $\End(V^{\ot 3})$.  The quantity $\theta_\mu:=q^{\lan\mu+2\rho,\mu\ra}$ is the \emph{twist} corresponding to the object $W$.  It follows from the axioms of a ribbon category that the scalar by which $(c_{V,V})^2$ acts on $\Hom(W,V^{\ot 2})$ is $\theta_W/(\theta_V)^2$ for $W$ a simple subobject of $V\ot V$ (with $V$ simple and $V\ot V$ multiplicity free). For a general ribbon category the objects have no structure so the determining the signs of the eigenvalues of $c_{V,V}$ is a bit more delicate.

\section{Weakly integrality}\label{3:wiclass}
In this section we classify the pairs $(\g,\ell)$ with $\g$ of exceptional Lie type such that $\CC(\g,q,\ell)$ is weakly integral.  The classification is summarized in Table \ref{weak}.

One important fact is the following:
\begin{prop}(see \cite{ENO})
If $\CC$ is a weakly integral fusion category then $\CC_{ad}$ is an integral fusion category.
\end{prop}
Here $\CC_{ad}$ is the \emph{adjoint subcategory} generated by simple subobjects of $X\ot X^*$ for each $X$.  Thus to show that $\CC$ is not weakly integral it is enough to find an object in $\CC_{ad}$ that has non-integer $FP$-dimension.

\begin{prop} The ribbon category $\CC(\g,q,\ell)$ is weakly integral with rank at least $2$ if and only if $(\g,\ell)$ is in Table \ref{weak}.
\end{prop}
\begin{table}\caption{Non-trivial Weakly Integral $\mathcal{C}(\g,q,\ell)$}\label{weak}
\begin{tabular}{*{3}{|c}|}
\hline
$\g$ & $\ell$ & Notes\\
\hline\hline
 $\e_6$& $13$ & pointed, rank $3$
\\ \hline  $\e_7$& $19$ & pointed, rank $2$
\\ \hline  $\e_8$& $32$ & rank $3$
\\ \hline  $\f_4$& none &
\\ \hline  $\g_2$& $8$ & pointed, rank $2$
\\
\hline
\end{tabular}
\end{table}

\begin{proof} Direct computation show that the categories listed are weakly integral.


To verify that no other weakly integral cases exist is tedious but straightforward.  We will first illustrate how this is done by giving complete details for type $\g_2$, and then describe how to carry out the computation for other types.

First suppose that $3\mid\ell$.  $\CC:=\CC(\g_2,q,\ell)$ is trivial (rank $1$) for $\ell=12$ so we assume $\ell\geq 15$.  Let $\la_1$ be the highest weight corresponding to the $7$-dimensional representation of $\g_2$, and denote by $X$ the corresponding simple object in $\CC$.  Note that $X$ is self-dual and $X\in\CC_{ad}$ since $X$ is a subobject of $X\ot X$.  Now $\FPdim(X)=\frac{[2][7][12]}{[4][6]}$, so if $\FPdim(X)=k\in\Z$ then the primitive $2\ell$th root of unity $q$ satisfies the polynomial in $\Z[q]$:
$$q^{20}+q^{18}+q^{12}+(1-k)q^{10}+q^8+q^2+1.$$
Thus the minimal polynomial of $q$ must divide the above polynomial so in particular, $\phi(2\ell)\leq 20$ where $\phi$ is Euler's totient function.  This implies that $\ell\in\{15,18,\ldots,33\}$.  Computing $\FPdim(X)$ for these 7 values of $\ell$ we find no integers, hence $\CC$ is not integral.

Now suppose that $3\nmid\ell$.  Here $\CC$ is trivial for $\ell=7$ and pointed (integral) for $\ell=8$ so we assume $\ell\geq 10$. From $\FPdim(X)=[7]$ one finds that if $\FPdim(X)=k\in\Z$ then $\phi(2\ell)\leq 12$ as $q$ must satisfy a degree $12$ polynomial.  From this we reduce the problem to checking that $\FPdim(X)$ is non-integral for $\ell\in\{10,11,13,14\}$, which is easily done.  Thus $\CC$ is not integral except in the case $\ell=8$ found in Table \ref{weak}.

The other exceptional types can be done in a similar fashion: find a convenient object $X\in\CC(\g,q,\ell)_{ad}$, assume that $\FPdim(X)=k\in\Z$ to obtain an upper bound on $\ell$ by bounding $\phi(2\ell)$ as above and then checking that $\FPdim(X)$ is non-integral for the (finitely many) remaining values of $\ell$.  Table \ref{nonint} gives the necessary data for all exceptional Lie types.  The second column of Table \ref{nonint} gives a pair of weights $(\nu,\mu)$ such that $V_{\nu}\subset V_{\mu}\ot V^*_{\mu}$, that is, so that $V_{\nu}\in\CC(\g,q,\ell)_{ad}$.  The corresponding $\FPdim(V_{\nu})$ is given in the third column, and the upper bound on $\ell$ that the assumption $\FPdim(V_\nu)\in\Z$ induces is found in the last column.

\end{proof}

\begin{table}\caption{Objects of Non-integral $FP$-dimension}\label{nonint}
\begin{tabular}{*{4}{|c}|}
\hline
$\g$ &  $(\nu,\mu)$: $V_{\nu}\subset V_{\mu}\ot V^*_{\mu}$ & $\FPdim(V_\nu)$ & maximum $\ell$\\
\hline\hline
 $\e_6$& $(\la_2,\la_1)$ & $\frac{[8][9][13]}{[4][3]}$ & $75$
\\ \hline  $\e_7$& $(\la_1,\la_7)$   & $\frac{[12][14][19]}{[4][6]}$ & $120$
\\ \hline  $\e_8$& $(\la_8,\la_8)$&  $\frac{[20][24][31]}{[6][10]}$ & $210$
\\ \hline  $\f_4$, $\ell$ even & $(\la_1,\la_1)$ & $\frac{[3][8][13][18]}{[4][6][9]}$ & $66$
\\ \hline  $\f_4$, $\ell$ odd & $(\la_1,\la_1)$ & $\frac{[13][8]}{[4]}$ & $51$
\\ \hline  $\g_2$, $3\mid\ell$ & $(\la_1,\la_1)$  & $\frac{[2][7][12]}{[4][6]}$ & $33$
\\ \hline  $\g_2$, $3\nmid\ell$ & $(\la_1,\la_1)$  &  $[7]$ & $14$
\\
\hline
\end{tabular}
\end{table}

\section{Failure of property $\mathbf{F}$}\label{4:noF}

In this section we demonstrate that property $\mathbf{F}$ fails for each pair $(\g,\ell)$ such that $\CC(\g,q,\ell)$ is not weakly integral (see Table \ref{weak}).

\subsection{Lie type $G_2$}

\begin{thm}
The non-trivial categories $\CC(\g_2,q,\ell)$ do not have property $\mathbf{F}$ unless $\ell=8$.
\end{thm}
\begin{proof}
We first remark that it is enough to consider only the specific choice $q=e^{\pi i/\ell}$ as Galois conjugation
does not affect the question of finiteness of the image of $\B_3$: the relations in a finite group presentation induce polynomial equations in the entries with integer coefficients.  By non-trivial we mean that the rank is at least $2$.

Let $V$ be the simple object labeled by $\la_1=\ve_1-\ve_3$ i.e. the highest weight of the $7$-dimensional fundamental representation of $\g_2$.  There are two cases to consider: $3\mid\ell$ and $3\nmid\ell$.

We first consider the cases where $3\mid \ell$.  Here $\ell=12$ corresponds to the trivial rank $1$ category and $\ell=15$ has rank $2$.  The latter case was considered in \cite{RSW}, where it is called the ``Fibonacci category" and can be identified with a subcategory of $\CC(\ssl_2,q,5)$.  In particular it is known to have infinite braid group image (going back to \cite{jones86}).  We claim that, provided $\ell\geq 18$, the rank of $\CC(\g_2,q,\ell)$ is at least $4$ and $\Hom(V,V^{\ot 3})$ is a $4$-dimensional irreducible representation of $\B_3$.  Note that $V^{\ot 2}\cong \one\oplus V\oplus V_{\la_2}\oplus V_{2\la_1}$, and that each of these summands is in $C_\ell(\g_2)$.  One computes the eigenvalues of the image of $\sigma_1$ acting on $\Hom(V,V^{\ot 3})$ using Reshetikhin's formula (\ref{eigs}): $\{1,-q^6,-q^{12},q^{14}\}$ (note that $\bigwedge^2 V\cong V\oplus V_{\la_2}$ which accounts for the signs).  It is convenient to rescale these to:  $S:=\{q^{-12},-q^{-6},-1,q^2\}$ which of course does not affect the projective image.  These eigenvalues are distinct for $q=e^{\pi i/\ell}$ with $\ell\geq 18$ and so it follows from Lemma \ref{twlemma} that $\B_3$ acts irreducibly on $\Hom(V,V^{\ot 3})$ and the projective order of the image of $\sigma_1$ is $2\ell$ if $\ell$ is odd and $\ell$ if $\ell$ is even.  Next we claim that the image of $\B_3$ is primitive: by Theorem \ref{rttheorem}(b) is is enough to checks that $S$ is not of the form $\{\pm r,\pm s\}$ or $\chi\{1,\omega,\omega^2\}\cup\{\alpha\}$
with $\omega$ a 3rd root of unity.  Since $q=e^{\pi i/\ell}$ with $\ell\geq 18$ we find that $1\not\in S$ so $S\not=\{\pm r,\pm s\}$.  For $\ell\geq 18$ none of $-1/q^2,q^{-12}/q^2$ or $-q^{-6}/q^2$ is equal to $1$ so we need only check that $\{-1,q^{-12},-q^{-6}\}\not=\chi\{1,\omega,\omega^2\}$ for any $\chi$.  Assuming to the contrary we must have $q^{-18}=1$ but this implies $q^2$ is a $9$th root of unity which contradicts $\ell\geq 18$.  Now it follows from Theorem \ref{rttheorem}(a) and (d)(iii)] that the image of $\B_3$ is infinite unless $\ell=24$.  For the case $\ell=24$ one must work a little harder: in this case there are, up to equivalence, two irreducible $4$ dimensional representations of $\B_3$ with these eigenvalues, explicitly described in \cite[Prop. 2.6]{TW}.  For example for one of the two choices the image of $\sigma_1$ is
$$A:=\begin{pmatrix} {q}^{-12}&{\frac {{q}^{8}+ {q}^{4}+1}{{q}^{6}}}&-{\frac {{q}^{8}+{q}^{4}+1}{{q}^{14}}}&-1
\\0&{q}^{2}&-{\frac {{q}^{4}-1}{{q}^{10}}}&-1\\0&0&-{q}^{-6}&-1
\\0&0&0&-1\end{pmatrix}$$ while
the image of $\sigma_2$ is:
$$B:=\begin{pmatrix}-1&0&0&0\\{q}^{-6}&-{q}^{-6}&0&0\\ {q}^{6}
&- \left( {q}^{4}+1 \right) {q}^{2}&
{q}^{2}&0\\-1 &{
\frac {{q}^{8} +{q}^{4}+1}{{q}^{8}}}&-{\frac {{q}^{8}+{q}^{4}+1}{{q}^{12}}}&{q}^{-12}\end{pmatrix}$$
Substituting $q=e^{\pi i/24}$ one finds that the matrix $C:=AB^{-1}$ has infinite order.  Indeed by \cite[Lemma 5.1]{RT} one need only verify that $C^j$ is not proportional to $I$ for $1\leq j\leq 24$.

Next we consider the case $3\nmid \ell$.  In this case $\CC(\g_2,q,\ell)$ is trivial for $\ell=7$ and is pointed of rank $2$ for $\ell=8$.  For $\ell\geq 10$ we again find that $\Hom(V,V^{\ot 3})$ is a $4$-dimensional irreducible $\B_3$-representation as we have the same decomposition of $V^{\ot 2}$ and eigenvalues as above.  The case $\ell=10$ (rank $4$) has been studied: $\CC(\g_2,q,10)$ has the Fibonacci category as a (modular) subcategory (see \cite[Theorem 3.4]{RJPAA}) and hence has infinite braid group image.  For $\ell\geq 11$, Theorem \ref{rttheorem}(d) is again sufficient to conclude the image of $\B_3$ is infinite except for the case $\ell=20$.  We may use the same matrices $A,B$ and $C$ described above and explicitly check that $C$ has infinite order for $q=e^{\pi i/20}$.
\end{proof}
\subsection{Lie type $F_4$}
\begin{thm}
The non-trivial categories $\CC(\mathfrak{f}_4,q,\ell)$ do not have property $\mathbf{F}$.
\end{thm}

\begin{proof}

Now let $\g=\mathfrak{f}_4$, and $V$ be the simple object in $\CC(\mathfrak{f}_4,q,\ell)$ analogous to the ($26$-dimensional) vector representation of $\mathfrak{f}_4$ labeled by $\la_1$.  Again there are two cases $2\mid\ell$ and $2\nmid\ell$.

First, suppose $\ell$ is even. Then if $22\leq\ell$, $\dim\Hom(V,V^{\ot 3})=5$ and the eigenvalues of the image of $\sigma_1$ are (up to an overall scale factor): $\{q^{-24},q^{-12},q^2,-1,-q^{-6}\}$.  Notice that these eigenvalues are distinct unless $\ell=24$ so that $\Hom(V,V^{\ot 3})$ is an irreducible $\B_3$-representation if $\ell\not=24$.  Moreover, the image is primitive (see Theorem \ref{rttheorem}(b)) so Theorem \ref{rttheorem}(d)(iv) implies that the image of $\B_3$ is infinite for $\ell=22$ or $26\leq\ell$ as the projective order of the image of $\sigma_1$ is $\ell$ in these cases. In the case $\ell=24$ we have repeated eigenvalues, hence either the representation is reducible (or has infinite image by Theorem \ref{rttheorem}(a)).  Thus we must consider the $\ell=24$ case separately.  In this case $\CC(\mathfrak{f}_4,q,24)$ has rank $9$: the simple objects are labeled by $$\{0,\la_1,2\la_1,3\la_1,\la_2,\la_3,\la_4,\la_1+\la_2,\la_1+\la_4\}.$$  If we let $U$ be the simple object labeled by $\la_4$ we compute that $U^{\ot 2}$ decomposes as $\one\oplus V_{2\la_1}\oplus U\oplus V_{\la_3}$.  To see this one compute the second tensor power of the fundamental $52$-dimensional representation with highest weight $\la_4$ and then discards the object labeled by $2\la_4$ as this lies on the upper hyperplane of the Weyl alcove for $\ell=24$.  That is, $\lan 2\la_4+\rho,\vartheta_0\ra=24$.  One Then computes the eigenvalues of the action of $\sigma_1$ on $\Hom(U,U^{\ot 3})$ to be: $\{1,q^{26},-q^{18},-q^{36}\}$.  For $q=e^{\pi i/24}$ these are distinct so that $\Hom(U,U^{\ot 3})$ is an irreducible $4$-dimensional $\B_3$-representation.  However, the projective order of the image of $\sigma_1$ is $24$ so we must resort to explicit computations as in the $\g_2$ situation above.  We again find that the image of $\sigma_1\sigma_2^{-1}$ has infinite image.

Next assume that $\ell$ is odd.  We have $\dim\Hom(V,V^{\ot 3})=5$ when $15\leq\ell$, and the projective order of the image of $\sigma_1$ is $2\ell$ so the image of $\B_3$ is again always infinite by Theorem \ref{rttheorem}(d).

\end{proof}

\subsection{Lie types $E_6,E_7$ and $E_8$}
\begin{thm}
 The non-trivial categories $\CC(\mathfrak{e}_N,q,\ell)$ do not have property $\mathbf{F}$ unless $(N,\ell)\in\{(6,13),(7,19),(8,32)\}$.
\end{thm}

\begin{proof}

The braid group representations for $U_q\mathfrak{e}_N$ with $q$ generic have been studied at length in \cite{wenzlEN}, and the results there can be used to give a uniform proof of the infinitude of the braid group image for $N=6,7$ and $8$.  In this case we take $V$ to be the simple object analogous to the vector representation of $\mathfrak{e}_N$.  The highest weight of $V$ corresponds to the node in the Dynkin diagram of $E_N$ furthest from the triple point (for $N=6$ this is ambiguous but we may pick either as they correspond to dual objects).  We will call this highest weight $\la_1$ (although this does not coincide with \cite{Bou} in some cases).

The $\B_3$ representation space we consider is $\Hom(V_{\la_1+\la_N},V^{\ot 3})$.  Note that $V_{\la_1+\la_N}$ appears in what Wenzl calls $V^{\ot 3}_{new}$ in \cite{wenzlEN}.
This space is $3$-dimensional provided $\la_1+\la_N\in C_\ell(\mathfrak{e}_N)$.  For $N=6$ this is satisfied for $\ell\geq 14$, for $N=7$ we need $\ell\geq 21$ and for $N=8$ we must have $\ell\geq 34$.  Let us say that $\ell$ is in the \textit{stable range} if $\ell$ is large enough to ensure $\la_1+\la_N\in C_\ell(\mathfrak{e}_N)$.  We will consider the cases corresponding to the pairs $(N,\ell)\in\{ (7,20),(8,33)\}$ separately.  Provided $\ell$ is in the stable range the braiding eigenvalues on $H:=\Hom(V_{\la_1+\la_N},V^{\ot 3})$ are (up to an overall scale factor): $\{q,-q^{-1},q^{3-2N}\}$ (see \cite{wenzlEN}) where $N=6,7$ or $8$.  In fact, it is observed in \cite[Remark 5.10]{wenzlEN} that this $3$-dimensional representation $H$ is equivalent to the $3$-dimensional representation of $\B_3$ obtained from BMW-algebras (see \cite{wenzlbcd,tubawenzlbcd}) specialized at $r=q^{2N-3}$.  The irreducibility of this representation at roots of unity is analyzed in \cite{wenzlbcd}.  In particular \cite[Theorem 6.4(a)]{wenzlbcd} implies that as long as $2N-3<\ell-2$ the corresponding $3$-dimensional $\B_3$ representation is irreducible.  This is clearly satisfied for $\ell$ in the stable range.  An alternative approach to showing that $H$ is irreducible is to follow the proof in \cite[Section 4]{wenzlEN}, which is a matter of checking that the argument there is valid for the root of unity case provided $\ell$ is in the stable range.

Having shown that $H$ is irreducible it is now a routine application of Theorem \ref{rttheorem}(c)(i) and (d)(ii) to see that the image of $\B_3$ on $H$ is primitive and hence infinite as the projective order of the image of $\sigma_1$ is at least $8$.  This proves the statement for all but the two cases $(N,\ell)\in\{(7,20),(8,33)\}$.

One finds that $\CC(\mathfrak{e}_7,q,20)$ has rank $6$ and is a product of two well-known categories: the Fibonacci category (rank $2$) and the Ising category (rank $3$).  The former has infinite braid group image (although the latter does not!) so we may deduce the infinitude of the braid group image for $(N,\ell)=(7,20)$.

The rank $5$ category $\CC(\mathfrak{e}_8,q,33)$ is conjugate to $\CC(\mathfrak{f}_4,q,22)$ and can also be realized as a modular subcategory of $\CC(\ssl_2,q,11)$.  We have already seen that $\CC(\mathfrak{f}_4,q,22)$ has infinite braid group image, so we have shown that the case $(N,\ell)=(8,33)$ has infinite braid group image finishing the proof.

\end{proof}

\section{Conclusions}\label{5:conclusions}
With these results the verification of the ``non-weakly-integral implies no property $\mathbf{F}$'' direction of Conjecture \ref{propfconj} for quantum group categories is essentially complete.  Although a classification of weakly integral categories of classical Lie type has not appeared (to our knowledge), it is essentially known to experts.  In type $A$ it is known that $\CC(\ssl_N,q,\ell)$ is weakly integral for $\ell\in\{N,N+1,4,6\}$ and for type $C$ we have $\CC(\mathfrak{sp}_4,q,10)$ weakly integral.  There are two infinite families of weakly integral categories coming from types $B$ and $D$: $\CC(\so_N,q,2N)$ with $N$ odd and $\CC(\so_M,q,M)$ with $M$ even (see \cite{NR}).  In addition $\CC(\so_N,q,2N-2)$ with $N$ odd and $\CC(\so_M,q,M-1)$ with $M$ even are weakly integral but are always rank $3$ and $4$ respectively and only give rise to finitely many inequivalent categories.  This is expected to be a complete list (except possibly for some low-rank coincidences).  The results of \cite{jones86,FLW,LRW} show that in all but these cases property $\mathbf{F}$ fails for $\CC(\g,q,\ell)$ of classical Lie type.  

We conclude with two ``plausibility arguments'' for the general conjecture.
Firstly, if $G$ is a primitive linear group of degree $m$ then the projective order of any element of $G$ is bounded by a ($G$-independent) function of $m$.  Indeed, see \cite[Corollary 4.3]{LRW} for a result of this type.  Now it is known that for a (unitary) modular category, the $\FPdim(X_j)$ lie in the cyclotomic field $\Q(\theta_1,\ldots,\theta_k)$ generated by the twists (see \cite{NS}).  So if $\FPdim(X_j)$ is far from being integral, that is $[\Q(\FPdim(X_j):\Q]$ is large, then the order of some $\theta_i$ is large.  Since the eigenvalues of the image of $\sigma_i$ are square-roots of products of twists it follows that the projective order of the image of $\sigma_1$ in some $\B_3$ representation is large.  On the other hand the degrees of the irreducible representations of $\B_3$ associated to an object $X$ are bounded by $\FPdim(X)^2$ (the maximum number of simple subobjects of $X^{\ot 2}$) and is typically much smaller.  Thus it is reasonable to expect that non-weakly-integral categories give rise to infinite braid group images.

The main result of \cite{ERW} is that group-theoretical braided fusion categories have property $\mathbf{F}$.  Here $\CC$ is \emph{group-theoretical} if the Drinfeld center of $\CC$ is equivalent to $\Rep(D^\om G)$ for some finite group $G$ and cocycle $\om$.  A group-theoretical fusion category is always integral, but not conversely.  This notion has been generalized to weakly group-theoretical fusion categories (see \cite{ENO2}).  Roughly speaking, weakly group-theoretical fusion categories are those that can be defined in terms of finite group data.  
Weakly group-theoretical fusion categories are always weakly integral, and no counterexample to the converse is currently known.  If indeed weak integrality and weak group-theoreticity are equivalent notions, one need only show that weakly group-theoretical braided fusion categories have property $\mathbf{F}$ to prove the other direction of Conjecture \ref{propfconj}.

\end{document}